\magnification 1200
\input plainenc
\input amssym
\fontencoding{T2A}
\inputencoding{utf-8}
\tolerance 4000
\relpenalty 10000
\binoppenalty 10000
\parindent 1.5em

\hsize 17truecm
\vsize 24truecm
\hoffset 0truecm
\voffset -0.5truecm

\font\TITLE labx1440
\font\tenrm larm1000
\font\cmtenrm cmr10
\font\tenit lati1000
\font\tenbf labx1000
\font\tentt latt1000
\font\teni cmmi10 \skewchar\teni '177
\font\tensy cmsy10 \skewchar\tensy '60
\font\tenex cmex10
\font\teneufm eufm10
\font\eightrm larm0800
\font\cmeightrm cmr8
\font\eightit lati0800
\font\eightbf labx0800
\font\eighttt latt0800
\font\eighti cmmi8 \skewchar\eighti '177
\font\eightsy cmsy8 \skewchar\eightsy '60
\font\eightex cmex8
\font\eighteufm eufm8

\font\cmsixrm cmr6

\font\sixbf labx0600
\font\sixi cmmi6 \skewchar\sixi '177
\font\sixsy cmsy6 \skewchar\sixsy '60
\font\sixeufm eufm6

\font\cmfiverm cmr5

\font\fivebf labx0500
\font\fivei cmmi5 \skewchar\fivei '177
\font\fivesy cmsy5 \skewchar\fivesy '60
\font\fiveeufm eufm5
\font\tencmmib cmmib10 \skewchar\tencmmib '177
\font\eightcmmib cmmib8 \skewchar\eightcmmib '177
\font\sevencmmib cmmib7 \skewchar\sevencmmib '177
\font\sixcmmib cmmib6 \skewchar\sixcmmib '177
\font\fivecmmib cmmib5 \skewchar\fivecmmib '177
\newfam\cmmibfam
\textfont\cmmibfam\tencmmib \scriptfont\cmmibfam\sevencmmib
\scriptscriptfont\cmmibfam\fivecmmib
\def\tenpoint{\def\rm{\fam0\tenrm}\def\it{\fam\itfam\tenit}%
	\def\bf{\fam\bffam\tenbf}\def\tt{\fam\ttfam\tentt}%
	\textfont0\cmtenrm \scriptfont0\cmsevenrm \scriptscriptfont0\cmfiverm
  	\textfont1\teni \scriptfont1\seveni \scriptscriptfont1\fivei
  	\textfont2\tensy \scriptfont2\sevensy \scriptscriptfont2\fivesy
  	\textfont3\tenex \scriptfont3\tenex \scriptscriptfont3\tenex
  	\textfont\itfam\tenit
	\textfont\bffam\tenbf \scriptfont\bffam\sevenbf
	\scriptscriptfont\bffam\fivebf
	\textfont\eufmfam\teneufm \scriptfont\eufmfam\seveneufm
	\scriptscriptfont\eufmfam\fiveeufm
	\textfont\cmmibfam\tencmmib \scriptfont\cmmibfam\sevencmmib
	\scriptscriptfont\cmmibfam\fivecmmib
	\normalbaselineskip 12pt
	\setbox\strutbox\hbox{\vrule height8.5pt depth3.5pt width0pt}%
	\normalbaselines\rm}
\def\eightpoint{\def\rm{\fam 0\eightrm}\def\it{\fam\itfam\eightit}%
	\def\bf{\fam\bffam\eightbf}\def\tt{\fam\ttfam\eighttt}%
	\textfont0\cmeightrm \scriptfont0\cmsixrm \scriptscriptfont0\cmfiverm
	\textfont1\eighti \scriptfont1\sixi \scriptscriptfont1\fivei
	\textfont2\eightsy \scriptfont2\sixsy \scriptscriptfont2\fivesy
	\textfont3\eightex \scriptfont3\eightex \scriptscriptfont3\eightex
	\textfont\itfam\eightit
	\textfont\bffam\eightbf \scriptfont\bffam\sixbf
	\scriptscriptfont\bffam\fivebf
	\textfont\eufmfam\eighteufm \scriptfont\eufmfam\sixeufm
	\scriptscriptfont\eufmfam\fiveeufm
	\textfont\cmmibfam\eightcmmib \scriptfont\cmmibfam\sixcmmib
	\scriptscriptfont\cmmibfam\fivecmmib
	\normalbaselineskip 11pt
	\abovedisplayskip 5pt
	\belowdisplayskip 5pt
	\setbox\strutbox\hbox{\vrule height7pt depth2pt width0pt}%
	\normalbaselines\rm
}

\def\No{\char 157}
\def\empty{}

\catcode`\@ 11
\catcode`\" 13
\def"#1{\ifx#1<\char 190\relax\else\ifx#1>\char 191\relax\else #1\fi\fi}

\def\newl@bel#1#2{\expandafter\def\csname l@#1\endcsname{#2}}
\openin 11\jobname .aux
\ifeof 11
	\closein 11\relax
\else
	\closein 11
	\input \jobname .aux
	\relax
\fi

\newcount\c@section
\newcount\c@subsection
\newcount\c@subsubsection
\newcount\c@equation
\newcount\c@bibl
\c@section 0
\c@subsection 0
\c@subsubsection 0
\c@equation 0
\c@bibl 0
\def\lab@l{}
\def\label#1{\immediate\write 11{\string\newl@bel{#1}{\lab@l}}%
	\ifhmode\unskip\fi}
\def\eqlabel#1{\rlap{$(\equation)$}\label{#1}}

\def\section#1{\global\advance\c@section 1
	{\par\vskip 3ex plus 0.5ex minus 0.1ex
	\rightskip 0pt plus 1fill\leftskip 0pt plus 1fill\noindent
	{\bf\S\thinspace\number\c@section .~#1}\par\penalty 25000%
	\vskip 1ex plus 0.25ex}
	\gdef\lab@l{\number\c@section.}
	\c@subsection 0
	\c@subsubsection 0
	\c@equation 0
}
\def\subsection{\global\advance\c@subsection 1
	\par\vskip 1ex plus 0.1ex minus 0.05ex{\bf\number\c@subsection. }%
	\gdef\lab@l{\number\c@section.\number\c@subsection}%
	\c@subsubsection 0\c@equation 0%
}
\def\subsubsection{\global\advance\c@subsubsection 1
	\par\vskip 1ex plus 0.1ex minus 0.05ex%
	{\bf\number\c@subsection.\number\c@subsubsection. }%
	\gdef\lab@l{\number\c@section.\number\c@subsection.%
		\number\c@subsubsection}%
}
\def\equation{\global\advance\c@equation 1
	\gdef\lab@l{\number\c@section.\number\c@subsection.%
	\number\c@equation}{\rm\number\c@equation}
}
\def\bibitem#1{\global\advance\c@bibl 1
	[\number\c@bibl]%
	\gdef\lab@l{\number\c@bibl}\label{#1}
}
\def\ref@ref#1.#2:{\def\REF@{#2}\ifx\REF@\empty{{\rm \S\thinspace#1}}%
	\else\ifnum #1=\c@section{{\rm #2}}\else{{\rm \S\thinspace#1.#2}}\fi\fi
}
\def\ref@eqref#1.#2.#3:{\ifnum #1=\c@section\ifnum #2=\c@subsection
	{{\rm (#3)}}\else{{\rm #2\thinspace(#3)}}\fi\else%
	{{\rm \S\thinspace#1.#2\thinspace(#3)}}\fi
}
\def\ref#1{\expandafter\ifx\csname l@#1\endcsname\relax
	{\bf ??}\else\edef\mur@{\csname l@#1\endcsname :}%
	{\expandafter\ref@ref\mur@}\fi
}
\def\eqref#1{\expandafter\ifx\csname l@#1\endcsname\relax
	{(\bf ??)}\else\edef\mur@{\csname l@#1\endcsname :}%
	{\expandafter\ref@eqref\mur@}\fi
}
\def\cite#1{\expandafter\ifx\csname l@#1\endcsname\relax
	{\bf ??}\else\hbox{\bf\csname l@#1\endcsname}\fi
}

\def\Wo{{\mathpalette\Wo@{}}W}
\def\Wo@#1{\setbox0\hbox{$#1 W$}\dimen@\ht0\dimen@ii\wd0\raise0.65\dimen@%
\rlap{\kern0.35\dimen@ii$#1{}^\circ$}}
\catcode`\" 12
\def\bolddelta{\mathchar"0\hexnumber@\cmmibfam 0E}
\catcode`\" 13
\catcode`\@ 12

\def\proof{\par\medskip{\rm Д$\,$о$\,$к$\,$а$\,$з$\,$а$\,$т$\,$е$\,$л$\,$ь%
	$\,$с$\,$т$\,$в$\,$о.}\ }
\def\endproof{{\parfillskip 0pt\hfill$\square$\par}\medskip}

\immediate\openout 11\jobname.aux


\frenchspacing\rm
\leftline{УДК~517.927+517.983.35}\vskip 0.25cm
{\leftskip 0cm plus 1fill\rightskip 0cm plus 1fill\parindent 0cm\baselineskip 15pt
\TITLE К вопросу о знакорегулярности положительных дифференциальных операторов
четвёртого порядка\par\vskip 0.25cm\rm А.$\,$А.~Владимиров%
\footnote{}{\eightpoint\rm Работа поддержана РНФ, грант \No$\,$14-11-00754.}\par}
\vskip 0.25cm
$$
	\vbox{\hsize 0.75\hsize\leftskip 0cm\rightskip 0cm
	\eightpoint\rm
	{\bf Аннотация:\/} Устанавливается, что положительность внутри открытого
	квадрата $(0,1)\times (0,1)$ функции Грина положительно определённого
	обыкно\-венного дифференциального оператора четвёртого порядка
	с распадающи\-мися гранич\-ными условиями представляет собой необходимое
	и достаточное условие того, чтобы этот оператор не понижал числа
	перемен знака.\par
	}
$$

\vskip 0.5cm
\section{Введение}\label{par:1}%
\subsection
Зафиксируем подпространство $\frak H\subseteq W_2^2[0,1]$, представляющее собой
ортогональное дополнение некоторого (возможно, пустого) линейно независимого набора
распределе\-ний класса $\mathop{\rm Lin}\{\bolddelta_0,\bolddelta_0'\}\cup
\mathop{\rm Lin}\{\bolddelta_1,\bolddelta_1'\}\subset W_2^{-2}[0,1]$. Рассмотрим
оператор $L\colon\frak H\to\frak H^*$, заданный тождеством
$$
	\langle Ly,z\rangle\equiv\int_0^1 p\,y''\overline{z''}\,dx+
		\langle q,\overline{y'}z'\rangle+\langle h,\overline yz\rangle,
	\leqno(\equation)
$$\label{eq:1.1}%
где функция $p\in L_\infty[0,1]$ равномерно положительна, а распределения
$q\in W_2^{-1}[0,1]$ и $h\in W_2^{-2}[0,1]$ вещественны (то есть сопоставляют
вещественнозначным пробным функциям вещественные значения). Целью настоящей статьи
является установление следующего факта.

\subsubsection\label{prop:1}
{\it Пусть оператор вида~\eqref{eq:1.1} положительно определён, а его функция Грина
$$
	G(t,s)\rightleftharpoons
		\overline{\langle\bolddelta_t,L^{-1}\bolddelta_s\rangle}
$$
положительна внутри открытого квадрата $(0,1)\times (0,1)$. Тогда число знакоперемен
никакой вещественнозначной функции $y\in\frak H$ не превосходит числа знакоперемен
соответствующего распределения $Ly\in\frak H^*$.
}

\bigskip
Здесь и далее мы считаем вещественную обобщённую функцию $f\in\frak H^*$
{\it имеющей не более $n\in\Bbb N$ перемен знака\/}, если она допускает в пространстве
$\frak H^*$ сколь угодно точную аппроксимацию имеющими не более $n$~перемен знака
непрерывными функциями.

Отметим, что функция Грина положительного оператора вида~\eqref{eq:1.1},
не понижающего числа перемен знака, заведомо является положительной внутри квадрата
$(0,1)\times (0,1)$. Тем самым, в рассматриваемом классе операторов
утверждение~\ref{prop:1} представляет собой критерий.

\subsection
Доказательство утверждения~\ref{prop:1} будет осуществлено нами на основе
последова\-тельной редукции задачи к следующей хорошо изученной ситуации.

\subsubsection\label{prop:1.1}
{\it Пусть оператор $L\colon W_2^2[0,1]\to W_2^{-2}[0,1]$ задан тождеством
$$
	\langle Ly,z\rangle\equiv\int_0^1 p\,y''\overline{z''}\,dx+
		\alpha\,y(0)\overline{z(0)}+\beta\,y'(0)\overline{z'(0)}+
		\gamma\,y'(1)\overline{z'(1)},
$$
где коэффициенты $\alpha$, $\beta$ и $\gamma$ положительны, а функция $p\in
L_\infty[0,1]$ равномерно положительна. Тогда этот оператор не понижает числа
перемен знака.
}

\bigskip
Изложению основных промежуточных шагов редукционного процесса будет
посвя\-щён~\ref{par:2}. В заключительном~\ref{par:3} будет завершено доказательство
утвер\-ждения~\ref{prop:1}, а также обсуждены некоторые примеры его применения.


\section{Вспомогательные утверждения}\label{par:2}%
\subsection
Для полноты изложения укажем сперва доказательство утверждения~\ref{prop:1.1}.
Заметим, что используемое нами определение числа знакоперемен обобщённой функции
позволяет ограничиться установлением того факта, что для всякой вещественнозначной
функции $f\in C[0,1]$ соответствующая функция $y\rightleftharpoons L^{-1}f$ имеет
не б\'{о}льшее число знакоперемен. Такая функция $y\in W_2^2[0,1]$, однако,
представляет собой (см., напри\-мер,~[\cite{Vl:2004}]) решение граничной задачи
$$
	\displaylines{(py'')''=f,\cr
	\hbox to \displaywidth{\eqlabel{eq:1.2}\hfill $(py'')'(0)+\alpha y(0)=0,$
		\hfill}\cr
	\hbox to \displaywidth{\eqlabel{eq:1.3}\hfill $(py'')(0)-\beta y'(0)=0,$
		\hfill}\cr
	\hbox to \displaywidth{\eqlabel{eq:1.4}\hfill $(py'')'(1)=0,$\hfill}\cr
	\hbox to \displaywidth{\eqlabel{eq:1.5}\hfill $(py'')(1)+\gamma y'(1)=0.$
		\hfill}}
$$
Рассмотрим набор точек $0<\xi_1<\ldots<\xi_{n+1}<1$, где $n\geqslant 1$,
удовлетворяющий при всех $k\in 1\,.\,.\,n$ неравенствам $y(\xi_k)y(\xi_{k+1})<0$.
Согласно теореме Лагранжа, для всех $k\in 1\,.\,.\,n$ найдутся точки
$\xi_{1,k}\in (\xi_k,\xi_{k+1})$ со свойством $y'(\xi_{1,k})y(\xi_{k+1})>0$.
При этом либо найдётся также точка $\xi_{1,0}\in (0,\xi_1)$ со свойством
$y'(\xi_{1,0})y(\xi_1)>0$, либо, ввиду граничного условия~\eqref{eq:1.2}, будет
выполняться неравенство $(py'')'(0)\cdot y(\xi_1)<0$. Положим
$m\rightleftharpoons 0$ в первом случае и $m\rightleftharpoons 1$~--- во втором.
Для всех $k\in m\,.\,.\,n-1$ найдутся точки $\xi_{2,k}\in (\xi_{1,k},
\xi_{1,k+1})$ со свойством $(py'')(\xi_{2,k})y'(\xi_{2,k+1})>0$. Аналогичным
образом, ввиду граничных условий~\eqref{eq:1.3} и~\eqref{eq:1.5}, найдутся точки
$\xi_{2,m-1}\in (0,\xi_{1,m})$ и $\xi_{2,n}\in (\xi_{1,n},1)$ со свойствами
$(py'')(\xi_{2,m-1})y'(\xi_{1,m})>0$ и $(py'')(\xi_{2,n})y'(\xi_{1,n})<0$. Далее,
для всех $k\in {m-1}\,.\,.\,n-1$ найдутся точки $\xi_{3,k}\in (\xi_{2,k},
\xi_{2,k+1})$ со свойством $(py'')'(\xi_{3,k})\cdot (py'')(\xi_{2,k+1})>0$.
При этом выполняется неравенство $(py'')(\xi_{2,0})y(\xi_1)<0$, что, ввиду
сказанного ранее, гарантирует в случае $m=1$ существование точки $\xi_{3,-1}\in
(0,\xi_{2,0})$ со свойством $(py'')'(\xi_{3,-1})\cdot (py'')(\xi_{2,0})>0$.
Наконец, для всех $k\in -1\,.\,.\,n-2$ найдутся точки $\xi_{4,k}\in (\xi_{3,k},
\xi_{3,k+1})$ со свойством $f(\xi_{4,k})\cdot (py'')'(\xi_{3,k+1})>0$.
Поскольку, ввиду граничного условия~\eqref{eq:1.4}, найдётся также точка
$\xi_{4,n-1}\in (\xi_{3,n-1},1)$ со свойством $f(\xi_{4,n-1})\cdot
(py'')'(\xi_{3,n-1})<0$, то функция $f$ имеет не менее $n$~перемен знака.
Тем самым, утверждение~\ref{prop:1.1} справедливо.

Равносильные утверждению~\ref{prop:1.1} положения установлены, в частности,
в разделе~[\cite{GK:1950}: Гл.~\hbox{III}, \S$\,$8].

\subsection
Вторым шагом нашего построения будет установление справедливости следующего факта.

\subsubsection\label{prop:2.1}
{\it Пусть положительный оператор $L\colon W_2^2[0,1]\to W_2^{-2}[0,1]$ задан
тождеством
$$
	\langle Ly,z\rangle\equiv\int_0^1 p\,y''\overline{z''}\,dx+
		\langle q,\overline{y'}z'\rangle+\alpha\,y(0)\overline{z(0)},
$$
где коэффициент $\alpha$ положителен, функция $p\in L_\infty[0,1]$ равномерно
положительна, а распределение $q\in W_2^{-1}[0,1]$ вещественно. Тогда этот оператор
не понижает числа перемен знака.
}%
\proof
В рассматриваемом случае является положительным оператор $S\colon W_2^1[0,1]\to
W_2^{-1}[0,1]$ вида
$$
	\langle Sy,z\rangle\equiv\int_0^1 p\,y'\overline{z'}\,dx+
		\langle q,\overline yz\rangle.
$$
Согласно теории Штурма (см., например, [\cite{Vl:2009}]), это означает наличие
постоянной $\omega>0$, для которой функция $\sigma\rightleftharpoons\omega
S^{-1}(\bolddelta_0+\bolddelta_1)\in W_2^1[0,1]$ равномерно положительна
и под\-чиняется равенству $\int_0^1\sigma\,dx=1$. Введём в рассмотрение оператор $V\colon
W_2^2[0,1]\to W_2^2[0,1]$ вида $Vy\rightleftharpoons y\circ\tau$, где положено
$\tau(x)\rightleftharpoons\int_0^x\sigma\,dt$. Для такого оператора легко
проверяется справедливость тождеств
$$
	\belowdisplayskip 2pt
	\leqalignno{[Vy]'&=(y'\circ\tau)\cdot\sigma,\cr
		[Vy]''&=(y''\circ\tau)\cdot\sigma^2+(y'\circ\tau)\cdot\sigma',
	}
$$
$$
	\abovedisplayskip 2pt
	\bigl|[Vy]''\bigr|^2=(|y''|^2\circ\tau)\cdot\sigma^4+\sigma'\cdot
		\left({\bigl|[Vy]'\bigr|^2\over\sigma}\right)'.
$$
Из них, в свою очередь, немедленно вытекают соотношения
$$
	\leqalignno{\langle V^*LVy,y\rangle&=\int_0^1\hat p\cdot |y''|^2\,dx+
		\langle S\sigma,{\bigl|[Vy]'\bigr|^2\over\sigma}\rangle+
		\alpha\cdot|y(0)|^2\cr
		&=\int_0^1\hat p\cdot |y''|^2\,dx+\alpha\cdot|y(0)|^2+
		\omega\sigma(0)\cdot|y'(0)|^2+\omega\sigma(1)\cdot |y'(1)|^2,
	}
$$
в которых функция $\hat p\in L_\infty[0,1]$ определяется условием $\hat p\circ\tau=
p\sigma^3$. С учётом хорошо известного для полуторалинейных форм принципа
поляризации [\cite{Ka:1972}: Гл.~I, $(6.11)$] последнее означает справедливость
тождества
$$
	\langle V^*LVy,z\rangle\equiv\int_0^1\hat p\,y''\overline{z''}\,dx+
		\alpha\,y(0)\overline{z(0)}+\beta\,y'(0)\overline{z'(0)}+
		\gamma\,y'(1)\overline{z'(1)},
$$
где положено $\beta\rightleftharpoons\omega\sigma(0)>0$
и $\gamma\rightleftharpoons\omega\sigma(1)>0$.

Итак, оператор $V^*LV$ имеет вид, указанный при формулировке
утверждения~\ref{prop:1.1}, а потому не понижает числа перемен знака.
Однако операторы $V$ и $V^*$ оче\-видным образом сохраняют такое число.
Соответственно, оператор $L$ также не понижает числа перемен знака,
что и требовалось доказать.
\endproof

Метод, использованный нами для доказательства утверждения~\ref{prop:2.1}, ранее
применялся в работах~[\cite{BAV:2006}] и~[\cite{BAVSh:2013}]. Его основная идея,
однако, содержится уже в работе~[\cite{LN:1958}: Theorem~12.1].

\subsection
Третьим шагом нашего построения будет установление справедливости следующе\-го факта.

\subsubsection\label{prop:2.2}
{\it Пусть справедливо равенство $\frak H=W_2^2[0,1]$, а положительно определённый
оператор вида~\eqref{eq:1.1} имеет функцию Грина, равномерно положительную на
замкнутом квадрате $[0,1]\times [0,1]$. Тогда этот оператор не понижает числа
перемен знака.
}%
\proof
Сделанные предположения гарантируют равномерную положи\-тельность функции
$\sigma\rightleftharpoons L^{-1}\bolddelta_0\in W_2^2[0,1]$, а потому
и ограниченную обратимость связанного с ней оператора $V\colon W_2^2[0,1]\to
W_2^2[0,1]$ вида $Vy\rightleftharpoons\sigma y$. Далее, при помощи рассуждений,
аналогичных проведённым в ходе доказательства утверждения~\ref{prop:2.1},
устанавливается справедливость равенств
$$
	\leqalignno{\langle V^*LVy,y\rangle&=\int_0^1\hat p\cdot|y''|^2\,dx+
		\langle\hat q,|y'|^2\rangle+\langle L\sigma,\sigma\,|y|^2\rangle\cr
		&=\int_0^1\hat p\cdot|y''|^2\,dx+\langle\hat q,|y'|^2\rangle+
		\alpha\cdot|y(0)|^2,
	}
$$
где положено $\hat p\rightleftharpoons p\sigma^2$ и $\alpha\rightleftharpoons
\sigma(0)>0$, а распределение $\hat q\in W_2^{-1}[0,1]$ имеет вид
$$
	\langle\hat q,y\rangle\equiv\langle q,\sigma^2y\rangle+
		\int_0^1 2p\cdot[2(\sigma')^2-\sigma\sigma'']\,\overline y\,dx+
		\int_0^1 p\cdot(\sigma^2)'\,\overline{y'}\,dx.
$$
С учётом принципа поляризации это означает, что оператор $V^*LV$ имеет вид,
указанный при формулировке утверждения~\ref{prop:2.1}, а потому не понижает числа
перемен знака. Тогда, ввиду сохранения числа знакоперемен при действии операторов
$V$ и $V^*$, оператор $L$ также не понижает числа перемен знака, что и требовалось
доказать.
\endproof


\section{Завершение доказательства и обсуждение}\label{par:3}%
\subsection
Итак, пусть выполнены условия утверждения~\ref{prop:1}. Зафиксируем произвольное
значение $\varepsilon\in (0,1/2)$. Ввиду положительности квадратичной формы
оператора $L$ на подпространстве $\frak M\rightleftharpoons L^{-1}\mathop{\rm Lin}
\{\bolddelta_\varepsilon,\bolddelta'_\varepsilon,\bolddelta_{1-\varepsilon},
\bolddelta'_{1-\varepsilon}\}$,  всякий набор из четырёх величин $y(\varepsilon)$,
$y'(\varepsilon)$, $y(1-\varepsilon)$ и $y'(1-\varepsilon)$ однозначно определяет
некоторую функцию $y\in\frak M$. Положительность оператора $L$ означает также,
что поведение указанной функции слева от точки $\varepsilon$ полностью определяется
выбором пары величин $y(\varepsilon)$ и $y'(\varepsilon)$, а её поведение справа
от точки $1-\varepsilon$~--- выбором пары величин $y(1-\varepsilon)$
и $y'(1-\varepsilon)$.

Пусть теперь оператор $I\colon W_2^2[\varepsilon,1-\varepsilon]\to\frak H$ таков,
что функция $Iy\in\frak H$ всегда совпадает с исходной функцией
$y\in W_2^2[\varepsilon,1-\varepsilon]$ на отрезке $[\varepsilon,1-\varepsilon]$,
и с некоторой функцией класса $\frak M$~--- вне этого отрезка. Учитывая
стандартные представления
$$
	\displaylines{\langle q,y\rangle\equiv-\int_0^1 Q\,\overline{y'}\,dx+
		\alpha\,\overline{y(0)},\qquad Q\in L_2[0,1],\;\alpha\in\Bbb R,\cr
		\langle h,y\rangle\equiv\int_0^1 H\,\overline{y''}\,dx+
		\beta\,\overline{y(0)}+\gamma\,\overline{y'(0)},\qquad
		H\in L_2[0,1],\;\beta,\gamma\in\Bbb R,
	}
$$
функционалов $q\in W_2^{-1}[0,1]$ и $h\in W_2^{-2}[0,1]$, устанавливаем
справедливость тождества
$$
	\langle I^*LIy,y\rangle\equiv
		\int_\varepsilon^{1-\varepsilon} p\,|y''|^2\,dx+
		\langle\hat q,|y'|^2\rangle+\langle\hat h,|y|^2\rangle,
$$
где $\hat q\in W_2^{-1}[\varepsilon,1-\varepsilon]$
и $\hat h\in W_2^{-2}[\varepsilon,1-\varepsilon]$ суть некоторые вещественные
распределения. При этом функция Грина оператора $I^*LI$ представляет собой
ограничение функции Грина оператора $L$ на квадрат $[\varepsilon,1-\varepsilon]
\times[\varepsilon,1-\varepsilon]$, и потому на этом квадрате равномерно
положительна. С учётом утверждения~\ref{prop:2.2} это означает, что для всякой
имеющей не более $n$~перемен знака функции $f\in C[0,1]$, тождественно обращающейся
в нуль за пределами отрезка $[\varepsilon,1-\varepsilon]$, соответствующая функция
$y\rightleftharpoons L^{-1}f$ имеет на отрезке $[\varepsilon,1-\varepsilon]$
не более $n$~знакоперемен. Произвольность выбора значения $\varepsilon\in (0,1/2)$
означает теперь, что оператор $L^{-1}$ не повышает числа перемен знака никакой
финитной функции $f\in C[0,1]$. Для завершения доказательства
утверждения~\ref{prop:1} остаётся лишь учесть непрерывный характер вложения
$L_1[0,1]\hookrightarrow\frak H^*$.

\subsection
Использованное нами определение~\eqref{eq:1.1} обыкновенного дифференциального
оператора четвёртого порядка с распадающимися граничными условиями было введено
в работе~[\cite{Vl:2004}]. Оно представляет собой естественный результат дальнейшего
развития предложенного в работе~[\cite{NZSh:1999}] подхода к определению оператора
Штурма--Лиувилля с коэффициентами-распределениями.

Отметим, что класс задач с операторами вида~\eqref{eq:1.1} включает целый ряд
поста\-новок, в традиционной записи имеющих вид многоточечных. Так, формально
трёх\-точечная задача вида
$$
	\displaylines{y^{(4)}+y=f,\cr y(0)=y'(0)=y''(1)+y'(1)+y(1)=
		y'''(1)-y'(1)-y(1)=0,\cr
	y(1/2+0)-y(1/2-0)=y'(1/2+0)-y'(1/2-0)=y''(1/2+0)-y''(1/2-0)={}\quad\cr
		\kern 3truecm {}=y'''(1/2+0)-y'''(1/2-0)+y(1/2)=0
	}
$$
в действительности представляет собой задачу на пространстве
$$
	\frak H=\{y\in W_2^2[0,1]\;:\;y(0)=y'(0)=0\}
$$
с оператором вида~\eqref{eq:1.1}, где $p\equiv 1$, $q=\bolddelta_1$
и $h=1+\bolddelta_{1/2}+\bolddelta_1-\bolddelta_1'$.

\subsection
Пусть оператор вида~\eqref{eq:1.1} положительно определён и не понижает числа перемен
знака. Из теоремы Шёнберга [\cite{GK:1950}: Гл.~V, Теорема~$4'$] и обусловленных
положительной определённостью оператора $L$ неравенств
$$
	G\pmatrix{x_1&x_2&\ldots&x_n\cr x_1&x_2&\ldots&x_n}>0,\qquad
	0<x_1<x_2<\ldots<x_n<1,
$$
немедленно вытекает справедливость при всяких $0<x_1<x_2<\ldots<x_n<1$ и $0<\xi_1<
\xi_2<\ldots<\xi_n<1$ неравенства
$$
	G\pmatrix{x_1&x_2&\ldots&x_n\cr \xi_1&\xi_2&\ldots&\xi_n}\geqslant 0.
$$
Для случая, когда распределения $q,h\in W_2^{-1}[0,1]$ представляют собой
неотрицательные линейные комбинации кусочно-непрерывных функций с дельта-функциями
Дирака, справедливость указанного неравенства была установлена в недавней
работе~[\cite{K:2015}].


\vskip 0.4cm
\eightpoint\rm
{\leftskip 0cm\rightskip 0cm plus 1fill\parindent 0cm
\bf Литература\par\penalty 20000}\vskip 0.25truecm\penalty 20000
\bibitem{Vl:2004} {\it А.$\,$А.~Владимиров\/}. О сходимости последовательностей
обыкновенных дифференциальных опера\-торов~// Матем.~заметки.~--- 2004.~--- Т.~75,
\No~6.~--- С.~941--943.

\bibitem{GK:1950} {\it Ф.$\,$Р.~Гантмахер, М.$\,$Г.~Крейн\/}. Осцилляционные
матрицы и ядра и малые колебания механи\-ческих систем, изд.~2.
М.-Л.: Гостехиздат, 1950.

\bibitem{Vl:2009} {\it А.$\,$А.~Владимиров\/}. К осцилляционной теории задачи
Штурма--Лиувилля с сингулярными коэф\-фициентами~// Журнал выч.~матем.
и матем.~физики.~--- 2009.~--- Т.~49, \No~9.~--- С.~1609--1621.

\bibitem{Ka:1972} {\it Т.~Като\/}. Теория возмущений линейных операторов.
М.: Мир, 1972.

\bibitem{BAV:2006} {\it Ж.~Бен Амара, А.$\,$А.~Владимиров\/}. Об осцилляции
собственных функций задачи четвёртого порядка со спектральным параметром
в граничном условии~// Фунд. и~прикл. матем.~--- 2006.~--- Т.~12, \No~4.~---
С.~41--52.

\bibitem{BAVSh:2013} {\it Ж.~Бен Амара, А.$\,$А.~Владимиров, А.$\,$А.~Шкаликов\/}.
Спектральные и осцилляционные свойства одного линейного пучка дифференциальных
операторов четвёртого порядка~// Матем.~заметки.~--- 2013.~--- Т.~94, \No~1.~---
С.~55--67.

\bibitem{LN:1958} {\it W.~Leighton, Z.~Nehari\/}. On the oscillation of solutions
of self-adjoint linear differential equations of the fourth order~// Trans.~of
AMS.~--- 1958.~--- V.~89.~--- P.~325--377.

\bibitem{NZSh:1999} {\it М.$\,$И.~Нейман-заде, А.$\,$А.~Шкаликов\/}.
Операторы Шрёдингера с сингулярными потенциалами из прост\-ранств
мультипликаторов~// Матем. заметки.~--- 1999.~--- Т.~66, \No~5.~--- С.~723--733.

\bibitem{K:2015} {\it Р.$\,$Ч.~Кулаев\/}. Об осцилляционности функции Грина
многоточечной краевой задачи для уравне\-ния четвёртого порядка~//
Дифф.~уравнения.~--- 2015.~--- Т.~51, \No~4.~--- С.~445--458.
\bye